\documentclass[12pt,reqno]{amsart}
\usepackage{amsmath,amssymb,geometry,color,tikz-cd}
\usepackage[backref,pagebackref,pdftex,hyperindex]{hyperref}
\usepackage[alphabetic]{amsrefs}
\usepackage{amssymb}% see geometry.pdf on how to lay out the page. There's lots.
\usepackage{bbm}
\usepackage{enumitem}
\usepackage{upgreek}

\newtheorem{proposition}{Proposition}[section]
\newtheorem{theorem}[proposition]{Theorem}
\newtheorem{lemma}[proposition]{Lemma}

\newtheorem{conjecture}[proposition]{Conjecture}

\theoremstyle{definition}
\newtheorem{remark}[proposition]{Remark}
\newtheorem{definition}[proposition]{Definition}

\title{On existence of bounded relative-global complements for Fano fibrations}
\author{Sung Rak Choi}
\author{Chuyu Zhou}

\address{Department of Mathematics, Yonsei University, 50 Yonsei-ro, Seodaemun-gu, Seoul 03722, Republic of Korea}
\email{sungrakc@yonsei.ac.kr}
\address{Department of Mathematics, Yonsei University, 50 Yonsei-ro, Seodaemun-gu, Seoul 03722, Republic of Korea}
\email{chuyuzhou1@gmail.com}

\date{} % delete this line to display the current date
\thanks{
		\newline
		\indent
The authors are partially supported by Samsung Science and Technology Foundation under Project Number SSTF-BA2302-03.		
		}

\newcommand{\lct}{{\rm {lct}}}

\newcommand{\bQ}{\mathbb{Q}}
\newcommand{\bR}{\mathbb{R}}

\newcommand{\bZ}{\mathbb{Z}}

\newcommand{\mO}{\mathcal{O}}

\newcommand{\mR}{\mathcal{R}}

\begin{document}

\begin{abstract}
For Fano fibrations with $\epsilon$-lc singularities of a fixed dimension, we show the existence of bounded relative-global complements. If the base of the fibration is of dimension one, we even show the existence of bounded relative-global klt complements.
\end{abstract}

\maketitle

\setcounter{tocdepth}{1}

\tableofcontents

\section{Introduction}\label{sec: intro}

This note concerns the relative-global complements of Fano fibrations. We work over an algebraically closed field of characteristic zero throughout the note. We freely refer to \cite{Kollar13, KM98} for various concepts of singularities in birational geometry, e.g. lc, klt, and $\epsilon$-lc singularities.

We say $f: X\to Z$ is a Fano fibration (resp. $\epsilon$-lc Fano fibration) if the following conditions are satisfied:
\begin{enumerate}
\item $X$ and $Z$ are both normal projective varieties and $f_*\mO_X=\mO_Z$,
\item $X$ is klt (resp. $\epsilon$-lc) and $-K_X$ is ample over $Z$.
\end{enumerate}
One can similarly define the concepts of Fano type fibrations and $\epsilon$-lc Fano type fibrations just by replacing the second condition with the following:  
\begin{enumerate}
\item[$(2)'$] there exists an effective $\bQ$-divisor $B$ on $X$ such that $-K_X-B$ is ample over $Z$ and $(X, B)$ admits klt (resp. $\epsilon$-lc) singularities.
\end{enumerate}

\medskip

Let $X\to Z$ be an $\epsilon$-lc Fano type fibration of dimension $d$ (i.e. $\dim\ X=d$). We say it admits a \textit{bounded relative-global complement} (resp. \textit{bounded relative-global klt complement}) if there exist a positive integer $N$ (depending only on $d,\epsilon$) and an effective $\bQ$-divisor $\Lambda$ on $X$ such that $N(K_X+\Lambda)\sim 0/Z$ and $(X, \Lambda)$ is lc (resp. klt).

We study the following question on the existence of bounded relative-global complements for $\epsilon$-lc Fano fibrations. This is a variant of  Shokurov's original conjecture on the existence of bounded local complements for Fano fibrations. See \cite[Conjecture 1.3]{Sho00} and \cite[Theorem 1.8]{Birkar19}.

\begin{conjecture}\label{conj: main}
Let $d$ be a positive integer and $\epsilon$ a positive real number. Then there exists a positive integer $N$ depending only on $d, \epsilon$ satisfying the following: If $f: X\to Z$ is an $\epsilon$-lc Fano fibration of dimension $d$, then there exists an effective $\bQ$-divisor $\Lambda$ on $X$ such that $N(K_X+\Lambda)\sim 0/Z$ and $(X, \Lambda)$ is klt.
\end{conjecture}

We list some known results towards the conjecture due to Birkar.

In \cite{Birkar19}, the following existence of bounded relative-local complement is obtained.

\begin{theorem}{\rm (\cite[Theorem 1.8]{Birkar19})}\label{thm: local1}
Let $d$ be a positive integer. Then there exists a positive integer $N$ depending only on $d$ satisfying the following:
If $X\to Z$ is a Fano fibration of dimension $d$, then for any closed point $z\in Z$, there is an effective $\bQ$-divisor $\Lambda$ on $X$ such that $N(K_X+\Lambda)\sim 0/z$ and $(X, \Lambda)$ is lc over $z$.
\end{theorem}

In the above theorem, $/z$ or over $z$ means over an open neighborhood of $z\in Z$.

In \cite{Birkar23}, a refinement of Theorem \ref{thm: local1} is obtained if the base is assumed to be a curve.

\begin{theorem}{\rm (\cite[Corollary 1.4]{Birkar23})}\label{thm: local2}
Let $d$ be a positive integer and $\epsilon$ a positive real number. Then there exists a positive integer $N$ depending only on $d, \epsilon$ satisfying the following.
Let $X\to Z$ be an $\epsilon$-lc Fano fibration of dimension $d$ and $Z$ is a curve. Then for any closed point $z\in Z$ there is an effective $\bQ$-divisor $\Lambda$ on $X$ such that $N(K_X+\Lambda)\sim 0/z$ and $(X, \Lambda)$ is klt over $z$.
\end{theorem}

Both Theorem \ref{thm: local1} and Theorem \ref{thm: local2} are about local complements, i.e. complements over a closed point, while Conjecture \ref{conj: main} concerns relative-global complements.

The main results of this note generalize Birkar's results above from local complements to relative-global complements  for $\epsilon$-lc Fano fibrations.

The following result answers Conjecture \ref{conj: main} in lc setting.

\begin{theorem}{\rm (Theorem \ref{thm: main})}\label{main: global lc}
Let $d$ be a positive integer and $\epsilon$ a positive real number. Then there exists a positive integer $N$ depending only on $d, \epsilon$ satisfying the following: Let $f: X\to Z$ be an $\epsilon$-lc Fano fibration of dimension $d$. Then there exists an effective $\bQ$-divisor $\Lambda$ on $X$ such that $N(K_X+\Lambda)\sim 0/Z$ and $(X, \Lambda)$ is lc.
\end{theorem}

The next result answers Conjecture \ref{conj: main} when the base is one dimensional.

\begin{theorem}{\rm (Theorem \ref{thm: curve})}\label{main: global klt}
Let $d$ be a positive integer and $\epsilon$ a positive real number. Then there exists a positive integer $N$ depending only on $d, \epsilon$ satisfying the following. Let $f: X\to Z$ be an $\epsilon$-lc Fano fibration of dimension $d$ and $Z$ is a curve. Then there exists an effective $\bQ$-divisor $\Lambda$ on $X$ such that $N(K_X+\Lambda)\sim 0/C$ and $(X, \Lambda)$ is klt.
\end{theorem}

We state Theorem \ref{main: global lc} and Theorem \ref{main: global klt} just for $\epsilon$-lc Fano fibrations. However, we will actually prove them for $\epsilon$-lc log Fano fibrations (see Section \ref{sec: curve}), which give a little more general results by allowing a boundary.

In the last section, we will prove Theorem \ref{main: global lc} for Fano fibrations (not just for $\epsilon$-lc Fano fibrations) by a more direct way (without using induction). However, as our main concerns are about bounded relative-global klt complements and the proof towards Theorem \ref{main: global lc} really reveals what ingredients we are missing to obtain such complements (see Remark \ref{rem: what we lack}), we do not present this more general result in Introduction. We refer the readers to the last section.

\noindent
\subsection*{Acknowledgement}

We would like to thank all the participants of the reading seminar on learning Birkar's work starting from December 2023 in Yonsei University. We especially thank Sungwook Jang and Donghyeon Kim for their enthusiasm for numerous presentations for the seminar.

\section{Global complement engine}\label{sec: 2}

In \cite{Birkar21}, Birkar constructs a complement which is controlled globally.

\begin{theorem}{\rm (\cite[Theorem 1.9]{Birkar21})}\label{thm: global engine}
Let $d$ be a positive integer and $\mR\subset [0,1]$ be a finite set of rational numbers. Then there exists a natural number $n$ depending only on $d, \mR$ satisfying the following. Assume
\begin{enumerate}
\item $(X, \Delta)$ is a projective lc pair of dimension $d$, where the coefficients of $\Delta$ are contained in $\mR$;
\item $M$ is a semi-ample Cartier divisor on $X$ defining a contraction $f: X\to Z$ such that $X$ is Fano type over $Z$;
\item $M-(K_X+\Delta)$ is nef and big, and 
\item $S$ is a non-klt center of $(X, \Delta)$ with $M|_S\equiv 0$.
\end{enumerate}
Then there is an effective $\bQ$-divisor $\Lambda\geq \Delta$ on $X$ such that $n(K_X+\Lambda)\sim (n+2)M$, and $(X, \Lambda)$ is lc over a neighborhood of $z:=f(S)$ (which is not necessarily a closed point). 
\end{theorem}

Birkar applies the above global complement engine to get some bounded relative-global complements in a special case, i.e. $(d,r,\epsilon)$-Fano type fibration in \cite[Definition 1.1]{Birkar22}. We recall the definition of $(d,r,\epsilon)$-Fano type fibration here for convenience.

\begin{definition}
Let $d,r$ be positive integers and $\epsilon$ be a positive real number. A $(d,r,\epsilon)$-Fano type fibration consists of a pair $(X, B)$ and a contraction $f: X\to Z$ such that we have the following:
\begin{enumerate}
\item $(X, B)$ is a projective $\epsilon$-lc pair of dimension $d$,
\item $K_X+B\sim_\bR f^*L$ for some $\bR$-divisor $L$,
\item $X$ is Fano type over $Z$, i.e. $-K_X$ is big over $Z$,
\item $A$ is a very ample line bundle on $Z$ with $A^{\dim Z}\leq r$, and
\item  $A-L$ is ample.
\end{enumerate}
\end{definition}

\begin{theorem}{\rm (\cite[Theorem 1.7]{Birkar22})}\label{thm: global special}
Let $d,r$ be positive integers, $\epsilon$ be a positive real number, and $\mR\subset [0,1]$ be a finite set of rational numbers. Then there exist positive numbers $m,n$ depending on $d,r,\epsilon, \mR$ satisfying the following: Assume $(X, B)\to Z$ is a $(d,r,\epsilon)$-Fano type fibration and that
\begin{enumerate}
\item we have $0\leq \Delta\leq B$ and the coefficients of $\Delta$ are in $\mR$, and
\item $-(K_X+\Delta)$ is big over $Z$.
\end{enumerate}
Then there is a $\bQ$-divisor $\Lambda\geq \Delta$ such that $n(K_X+\Lambda)\sim mf^*A$ and $(X, \Lambda)$ is klt.
\end{theorem}

We remark here that Birkar obtains bounded relative-global klt complements in Theorem \ref{thm: global special} for $(d,r,\epsilon)$-Fano type fibrations, while we have bounded relative-global (lc) complements in Theorem \ref{main: global lc}  for more general objects, i.e. $\epsilon$-lc Fano fibrations (and even Fano fibrations in the last section). The advantage of $(d,r,\epsilon)$-Fano type fibrations lies in the fact that it admits boundedness property (see \cite[Theorem 1.2]{Birkar22}) and good singularity property for $\bQ$-linear system $|B+f^*A|_\bQ$ (see \cite[Theorem 1.4]{Birkar22}), which provide a large room for perturbation to get the bounded klt complements.

\section{Curve base}\label{sec: curve}

In this section, we prove Conjecture \ref{conj: main} when the base is of dimension one. We introduce the following notation. If $Y\to X$ is a proper birational contraction and $B$ is a $\bQ$-divisor on $X$, then we denote by $\widetilde{B}$  the birational transformation of $B$ on $Y$; if $X\dashrightarrow X'$ (resp. $X\dashrightarrow X''$) is a birational contraction map and $B$ is a $\bQ$-divisor on $X$, then we denote by $B'$ (resp. $B''$) for the pushforward of $B$ on $X'$.

From now on, we will deal with log Fano fibrations for generality. 
We say $f: (X,\Delta)\to Z$ is a log Fano fibration (resp. $\epsilon$-lc log Fano fibration) if the following conditions are satisfied:
\begin{enumerate}
\item $X$ and $Z$ are both normal projective varieties and $f_*\mO_X=\mO_Z$,
\item $(X, \Delta)$ is klt (resp. $\epsilon$-lc) and $-K_X-\Delta$ is ample over $Z$.
\end{enumerate}

We first treat a special case, where all the fibers are irreducible.

\begin{lemma}\label{lem: curve special}
Let $d$ be a positive integer, $\epsilon$ be a positive real number, and $\mR\subset[0,1]$ be a finite set of rational numbers. Then there exists a positive integer $N$ depending only on $d, \epsilon, \mR$ satisfying the following:
Let $f: (X, \Delta)\to C$ be an $\epsilon$-lc log Fano fibration of dimension $d$ where the base $C$ is a curve and the coefficients of $\Delta$ are contained in $\mR$. Suppose that all the fibers of $f$ are irreducible. Then there exists an effective $\bQ$-divisor $\Lambda\geq \Delta$ on $X$ such that $N(K_X+\Lambda)\sim 0/C$ and $(X, \Lambda)$ is klt.
\end{lemma}

\begin{proof}
We divide the proof into several steps.

\

\textit{Step 1}. Fix a closed point $z\in C$. Write $t_z:=\lct(f^*z; X, \Delta)\leq 1$ and let $T$ be an lc place of $(X, \Delta+t_zf^*z)$. By \cite{BCHM10}, there exists a birational contraction $g: Y\to X$ which extracts $T$, and we have
$$K_Y+\widetilde{\Delta}+t_z\widetilde{f^*z}+T=g^*(K_X+\Delta+t_zf^*z). $$
We will just assume that $T$ is $g$-exceptional since otherwise we may let $g$ be the identity in the proof below.
By taking a small $\bQ$-factorial modification of $Y$, we may assume that $Y$ is $\bQ$-factorial. By \cite[Theorem 1.1]{Birkar23}, there exists a rational number $t_0>0$ depending only on $d, \epsilon$ such that $t_z> t_0$. Considering the log pair $(Y, \widetilde{\Delta}+t_0\widetilde{f^*z}+T)$, we see
$$-K_Y-\widetilde{\Delta}-t_0\widetilde{f^*z}-T=g^*(-K_X-\Delta-t_zf^*z)+(t_z-t_0)\widetilde{f^*z} $$
is big over $C$. Running MMP over $C$ on 
$$-K_Y-\widetilde{\Delta}-t_0\widetilde{f^*z}-T,$$ we get a minimal model $Y\dashrightarrow Y'/C$ and $$-K_{Y'}-\widetilde{\Delta}'-t_0(\widetilde{f^*z})'-T'$$ 
is big and nef over $C$. We claim that $T$ is not contracted by $Y\dashrightarrow Y'$. Let $$0<P\sim_\bQ -K_X-\Delta/C$$ 
be a general ample $\bQ$-divisor on $X$ such that $(X, \Delta+t_zf^*z+P)$ is still lc. Then we have
$$K_Y+\widetilde{\Delta}+t_z\widetilde{f^*z}+T+g^*P=g^*(K_X+\Delta+t_zf^*z+P) $$
and $(Y, \widetilde{\Delta}+t_z\widetilde{f^*z}+T+g^*P)$ is lc. Since 
$$K_X+\Delta+t_zf^*z+P\sim_\bQ 0/C,$$ 
it is not hard to see that 
$$(Y', \widetilde{\Delta}'+t_z(\widetilde{f^*z})'+T'+(g^*P)')$$ 
is lc and so is $(Y', \widetilde{\Delta}'+t_0(\widetilde{f^*z})'+T')$. Suppose $T$ is contracted, i.e. $T'=0$. Since $Y\dashrightarrow Y'$ is obtained by MMP/$C$ on $-K_Y-\Delta-t_0\widetilde{f^*z}-T$, we have
$$0=A_{(Y, \Delta+t_0\widetilde{f^*z}+T)}(T)>A_{(Y', \widetilde{\Delta}'+t_0(\widetilde{f^*z})')}(T). $$
This means that $(Y', \widetilde{\Delta}'+t_0(\widetilde{f^*z})')$ is not lc and this is a contradiction.

\

\textit{Step 2}. 
We will apply Theorem \ref{thm: global engine} to 
$$g': (Y', \widetilde{\Delta}'+t_0(\widetilde{f^*z})'+T')\to C.$$
Recall that $-K_{Y'}-\widetilde{\Delta}'-t_0(\widetilde{f^*z})'-T'$ is big and nef over $C$.
Let $A$ be a very ample Cartier divisor on $C$ such that 
$$g'^*A-K_{Y'}-\widetilde{\Delta}'-t_0(\widetilde{f^*z})'-T'$$ 
is big and nef. By Theorem \ref{thm: global engine}, there exist a positive integer $n$ (depending only on $d,\epsilon$, $\mR$) and an effective $\bQ$-divisor $\Lambda\geq \widetilde{\Delta}'+t_0(\widetilde{f^*z})'+T'$ such that 
$$n(K_{Y'}+\Lambda)\sim (n+2)g'^*A$$ 
and $(Y', \Lambda)$ is lc over $z$. By \cite[6.1(3)]{Birkar19}, there exists an effective $\bQ$-divisor $\Gamma\geq \widetilde{\Delta}+t_0\widetilde{f^*z}+T$ on $Y$ such that 
$$n(K_Y+\Gamma)\sim (n+2)g^*f^*A$$ 
and $(Y, \Gamma)$ is lc over $z$. Denote by 
$\Gamma_z:=g_*\Gamma-\Delta-t_0f^*z$, we see 
$$n(K_X+\Gamma_z+\Delta+t_0f^*z)\sim (n+2)f^*A,$$ 
and $(X, \Gamma_z+\Delta+t_0f^*z)$ is lc over $z$. Up to replacing $n$ and $n+2$ by adding a common multiple, we may assume $nt_0$ is an integer and $n\Delta$ is integral.

\

\textit{Step 3}. 
%Under the same notations as in Step 2, we first see 
%$$n\Gamma_z\in |(n+2)f^*A-n(K_X+\Delta)-nt_0f^*z|. $$
Since $(X, \Delta)\to C$ is an $\epsilon$-lc log Fano fibration, by \cite[Theorem 1.1]{Birkar21}, there exists a positive integer $m$ depending only on $d, \epsilon, \mR$ such that $-m(K_X+\Delta)$ is a very ample line bundle over an open subset of $C$. We claim that, after replacing $m$ with a bounded multiple if necessary, one can find a Cartier divisor $B$ on $C$ such that there is an element
$$\Theta\in |f^*B-m(K_X+\Delta)|$$
satisfying that $(X, \Delta+\frac{1}{m}\Theta)$ is klt over the generic point of $C$. To see this, first note that $-m(K_X+\Delta)$ is very ample over an open subset $U\subset C$ and all the fibers are irreducible. Now there is a Cartier divisor $D$ on $f^{-1}(U)$ such that $(f^{-1}U, \Delta|_{f^{-1}U}+\frac{1}{m}D)$ is klt and 
$$\bar{D}+m(K_X+\Delta)\sim f^*J$$ 
for some $\bQ$-divisor $J$ on $C$, where $\bar{D}$ is the closure of $D$. Since $\bar{D}+m(K_X+\Delta)$ is a Weil divisor on $X$, so is $f^*J$. By \cite[Corollary 1.3]{Birkar23}, the multiplicities of fibers of $X\to C$ are bounded above by a positive integer $l$ depending only on $d, \epsilon, \mR$. Thus we see that the coefficients of $f^*J$ are contained in $\frac{1}{l!}\bZ$. Therefore, we have 
$$l!\bar{D}\in |l!f^*J-l!m(K_X+\Delta)| .$$ 
Now by replacing $m$ by $l!m$, taking $B=l!J$ and $\Theta=l!\bar{D}$, the claim is proven.

We may assume that $m$ is a multiple of $n$ and write $m=rn$. We can also assume that $B-r(n+2)A+rnt_0z$ is very ample by taking $B$ very ample.
Recall from Step 2 that
$$n\Gamma_z\in |(n+2)f^*A-n(K_X+\Delta)-nt_0f^*z|. $$
We note the following
$$rn\Gamma_z+f^*B-r(n+2)f^*A+rnt_0f^*z\in |f^*B-m(K_X+\Delta)|. $$
Recall also from Step $2$ that $(X, \Gamma_z+\Delta+t_0f^*z)$ is lc over $z$. By replacing $\Gamma_z$ with $\widehat{\Gamma}_z$ where
$$\Theta=m\widehat{\Gamma}_z \in |f^*B-m(K_X+\Delta)|$$
is a general element, we obtain a pair $(X, \widehat{\Gamma}_z+\Delta+t_0f^*z)$ which is lc over $z$. It is also klt over the generic point of $C$. Hence $(X, \widehat{\Gamma}_z+\Delta)$ is klt over $z$ and $$m(K_X+\widehat{\Gamma}_z+\Delta)\sim f^*B.$$

\

\textit{Step 4}. By Step 3, we could choose finitely many closed points $z_1, z_2,...,z_k$ such that for each $z_i$, there exist a very ample Cartier divisor $B_i$ on $C$ and an effective $\bQ$-divisor $\Gamma_i$ on $X$ such that
\begin{enumerate}
\item $m(K_X+\Gamma_i+\Delta)\sim f^*B_i$, and
\item $(X, \Gamma_i+\Delta)$ is klt over a neighborhood $U_i$ of $z_i$.
\end{enumerate}
We may assume $\cup_i U_i=C$. Note that 
$$m\Gamma_i\in |f^*B_i-m(K_X+\Delta)|. $$
Put
$$ m\Lambda_i:=m\Gamma_i+\sum_{j\neq i} f^*B_j.$$
We see
$$m\Lambda_i\in |\sum_jf^*B_j-m(K_X+\Delta)|,$$
and $m(K_X+\Lambda_i+\Delta)\sim 0/ C$. Since $(X, \Gamma_i+\Delta)$ is klt over $U_i$ and $\sum_{j\ne i}B_j$ is very ample, replacing $m\Lambda_i$ by a general element in the linear system
$$m\bar{\Lambda}\in|\sum_jf^*B_j-m(K_X+\Delta)|,$$
we obtain a pair $(X, \bar{\Lambda}+\Delta)$ which is klt over $U_i$ for each $i$ and $$m(K_X+\bar{\Lambda}+\Delta)\sim 0/C.$$ 
Recalling that $\cup_i U_i=C$, it is clear that $\bar{\Lambda}+\Delta$ is the required bounded relative-global klt complement.
\end{proof}

Now we are ready to remove the irreducibility condition on fibers in Lemma \ref{lem: curve special} and prove the full version for curve base. 

\begin{theorem}\label{thm: curve}
Let $d$ be a positive integer, $\epsilon$ be a positive real number, and $\mR\subset[0,1]$ be a finite set of rational numbers. Then there exists a positive integer $N$ depending only on $d, \epsilon, \mR$ satisfying the following. 
Let $f: (X, \Delta)\to C$ be an $\epsilon$-lc log Fano fibration of dimension $d$ where the base is a curve and the coefficients of $\Delta$ are contained in $\mR$. Then there exists an effective $\bQ$-divisor $\Lambda\geq \Delta$ on $X$ satisfying that $N(K_X+\Lambda)\sim 0/C$ and $(X, \Lambda)$ is klt.
\end{theorem}

\begin{proof}
We divide the proof into several steps.

\

\textit{Step 1}. Fix a closed point $z\in C$. Applying Step 1 and Step 2 of the proof of Lemma \ref{lem: curve special}, we know there exist a positive integer $n$ (depending only on $d, \epsilon, \mR$), a positive rational number $t_0$ (depending only on $d, \epsilon,\mR$), and an effective $\bQ$-divisor $\Gamma_z$ on $X$ such that
$$n(K_X+\Gamma_z+\Delta+t_0f^*z)\sim (n+2)f^*A $$
for some very ample Cartier divisor $A$ on $C$, and $(X, \Gamma_z+\Delta+t_0f^*z)$ is lc over $z$. It is clear to see the following
$$n\Gamma_z\in |(n+2)f^*A-n(K_X+\Delta)-nt_0f^*z| .$$
Up to replacing $n$ and $n+2$ by adding a common multiple, we may assume $nt_0$ is an integer and $n\Delta$ is integral.

\

\textit{Step 2}. By \textit{Step 1}, we could choose finitely many closed points $z_1, z_2,...,z_k$ such that for each $z_i$, there exist a very ample Cartier divisor $A_i$ on $C$ and an effective $\bQ$-divisor $\Gamma_i$ on $X$ such that
\begin{enumerate}
\item $n\Gamma_i\in |(n+2)f^*A_i-n(K_X+\Delta)-nt_0f^*z_i| $, and
\item $(X, \Gamma_i+\Delta+t_0f^*z_i)$ is lc over a neighborhood $U_i$ of $z_i$.
\end{enumerate}
We may assume $\cup_i U_i=C$. Let $B$ be a very ample Cartier divisor on $C$ such that $B-(n+2)A_i+nt_0z_i$ is very ample for every $i$.  Then
$$n\Lambda_i:=n\Gamma_i+f^*B-(n+2)f^*A_i+nt_0f^*z_i \in |f^*B-n(K_X+\Delta)|.$$
Since $(X, \Gamma_i+\Delta+t_0f^*z_i)$ is lc over $U_i$ and $B-(n+2)A_i+nt_0z_i$ is a very ample Cartier divisor for every $i$,  replacing $n\Lambda_i$ with a general element 
$$n\bar{\Lambda}\in |f^*B-n(K_X+\Delta)|, $$
we still have a piar $(X, \bar{\Lambda}+\Delta+t_0f^*z_i)$ which is lc over $U_i$ for each $i$ and 
$$n(K_X+\bar{\Lambda}+\Delta)\sim 0/C.$$ 
Since $\cup_i U_i=C$, it is clear that $\bar{\Lambda}+\Delta$ is a bounded relative-global complement.

\

\textit{Step 3}. Notation as in \textit{Step 2}, it is not hard to see that $(X, \bar{\Lambda}+\Delta+t_0\sum_{i=1}^k f^*z_i)$ is lc. We may assume that all non-irreducible fibers of $X\to C$ are contained in $\{f^*z_1,...,f^*z_k\}$.
Up to a small $\bQ$-factorialization, we assume that $(X, \bar{\Lambda}+\Delta+t_0\sum_i f^*z_i)$ is a $\bQ$-factorial log pair. Suppose  that 
$$f^*z_j:=\sum_{r\geq 1} n_{jr}T_{jr}, j=1,...,s$$ are all non-irreducible fibers. Then we can run an MMP/$C$ on 
$$\sum_{j=1}^s\sum_{r\ne 1}n_{jr}T_{jr}$$ 
to get a minimal model $X\dashrightarrow X'/C$, and all fibers of $f': X'\to C$ are irreducible (by Zariski Lemma). In particular, $f'^*z_j=n_{j1}T_{j1}'$ for $j=1,...,s$. We continue to run MMP/$C$ on $-K_{X'}-\Delta'$ to get the ample model of $-K_{X'}-\Delta'$, denoted by $X'\dashrightarrow X''/C$, where $-K_{X''}-\Delta''$ is ample over $C$. Note that $f'': (X'',\Delta'')\to C$ is still an $\epsilon$-lc log Fano fibration. Applying Lemma \ref{lem: curve special} we know there exist a positive integer $m$ (depending only on $d, \epsilon,\mR$) and an effective $\bQ$-divisor $\hat{\Gamma}\geq \Delta''$ on $X''$ such that $m(K_{X''}+\hat{\Gamma})\sim 0/C$ and $(X'', \hat{\Gamma})$ is klt.

Let $W$ be the common resolution of $X, X''$ given by $p: W\to X$ and $q: W\to X''$. Put
$$K_X+\hat{\Lambda}:=p_*q^*(K_{X''}+\hat{\Gamma}).$$
Then we see that $m(K_X+\hat{\Lambda})\sim 0/C$, $\hat{\Lambda}\geq \Delta$, and $(X, \hat{\Lambda})$ is sub-klt, and the negative part of $\hat{\Lambda}$ is supported on $\sum_{j=1}^s\sum_{r\neq 1}T_{jr}$. 
By now, we have two divisors on $X$, i.e. $\bar{\Lambda}$ and $\hat{\Lambda}$, satisfying the following conditions:
\begin{enumerate}
\item $n(K_X+\bar{\Lambda}+\Delta)\sim f^*H_1$ for some Cartier divisor $H_1$ on $C$, and and $(X, \bar{\Lambda}+\Delta+t_0\sum_if^*z_i)$ is lc,
\item $m(K_X+\hat{\Lambda})\sim f^*H_2$ for some Cartier divisor $H_2$ on $C$, $\hat{\Lambda}\geq \Delta$, and $(X, \hat{\Lambda})$ is sub-klt.
\end{enumerate}

\

\textit{Step 4}. We may assume that  $m$ is a multiple of $n$ and write $m=en$ for some positive integer $e$. Let $H$ be a very ample Cartier divisor on $C$ such that $m\hat{\Lambda}+f^*H-f^*H_2$ is effective and $H-eH_1$ is a very ample Cartier divisor. Then we have
$$m\bar{\Lambda}+f^*H-f^*eH_1\in |f^*H-m(K_X+\Delta)|$$ 
and
$$m(\hat{\Lambda}-\Delta)+f^*H-f^*H_2\in |f^*H-m(K_X+\Delta)|. $$
Taking a general element 
$$m\Theta\in |f^*H-m(K_X+\Delta)|, $$
we see that 
$$m(K_X+\Theta+\Delta)\sim 0/C$$ 
and $(X, \Theta+\Delta)$ is globally lc while it is klt over the generic point of $C$.

Let $y_1,...,y_l$ be the finitely many closed points on $C$ over which $(X, \Theta+\Delta)$ loses klt singularities. Applying Step 1 and Step 2 again, we could find another effective $\bQ$-divisor $\Omega$ on $X$ satisfying 
$$n(K_X+\Omega+\Delta)\sim 0/C$$ 
and $(X, \Omega+\Delta+t_0f^*y_i)$ is lc for every $1\leq i\leq l$. Finally we easily get a bounded relative-global klt complement by combining $\Theta$ and $\Omega$, say $\frac{\Theta+\Omega}{2}+\Delta$, i.e. 
$$2m\left(K_X+\frac{\Theta+\Omega}{2}+\Delta\right)\sim 0/C$$ 
and $(X, \frac{\Theta+\Omega}{2}+\Delta)$ is klt.
\end{proof}

\section{Higher dimensional base}\label{sec: 4}

In this section, we prove Conjecture \ref{conj: main} in lc setting.

\begin{theorem}\label{thm: main}
Let $d$ be a positive integer, $\epsilon$ be a positive real number, and $\mR\subset[0,1]$ be a finite set of rational numbers. Then there exists a positive integer $N$ depending only on $d, \epsilon, \mR$ satisfying the following. 
Let $f: (X, \Delta)\to Z$ be an $\epsilon$-lc log Fano fibration of dimension $d$ where the coefficients of $\Delta$ are contained in $\mR$. Then there exists an effective $\bQ$-divisor $\Lambda\geq \Delta$ on $X$ satisfying that $N(K_X+\Lambda)\sim 0/Z$ and $(X, \Lambda)$ is lc.
\end{theorem}

\begin{proof}
We divide the proof into several steps.

\

\textit{Step 1}. We already proved the theorem when the base is one dimensional (see Theorem \ref{thm: curve}), and we aim to prove by induction on the dimension of the base. Let $A$ be a very ample Cartier divisor on $Z$, and $H\in |A|$ a general element. Write $G:=f^*H$ and we still denote by $f$ for the restriction $G\to H$. Then $(G, \Delta|_G)\to H$ is also an $\epsilon$-lc log Fano fibration by the generality of $H$. By induction, there exist a positive integer $n$ (depending only on $d, \epsilon,\mR$) and an effective $\bQ$-divisor $B_G$ on $G$ such that 
$$n(K_G+B_G+\Delta|_G)\sim 0/H$$ and $(G, B_G+\Delta|_G)$ is lc.
Let $D$ be a Cartier divisor on $H$ such that 
$$n(K_G+B_G+\Delta|_G)\sim f^*D.$$ 
Then for some very ample Cartier divisor $L$ on $Z$ such that $L|_H-D$ is also a very ample Cartier divisor on $H$, we have 
$$n\left(K_G+B_G+\Delta|_G+\frac{f^*(L|_H)-f^*D}{n}\right)\sim  f^*(L|_H)=(f^*L)|_G.$$
For a general element
$$nF_G\in |nB_G+f^*(L|_H)-f^*D|,$$
we have
$$n(K_G+F_G+\Delta|_G)\sim (f^*L)|_G.$$ 
Thus  we have 
$$nF_G\in |-n(K_G+\Delta|_G+(f^*L)|_G)|,$$
and $(G, F_G+\Delta|_G)$ is also lc.
Consider the following exact sequence:
\begin{align*}
0\to \mO_X(-n(K_X+\Delta+G)-G+f^*L)\to &\mO_X(-n(K_X+\Delta+G)+f^*L)\\
\to &\mO_G(-n(K_G+\Delta|_G)+(f^*L)|_G)\to 0.
\end{align*}
We may choose $L$ sufficiently ample such that 
$$f^*L-(n+1)G-n(K_X+\Delta)-(K_X+\Delta)=f^*(L-(n+1)A)-(n+1)(K_X+\Delta) $$
is ample. By Kawamata-Viehweg vanishing, we have
$$H^1(X, -n(K_X+\Delta+G)-G+f^*L)=0, $$ 
which leads to the following surjective map:
$$H^0(X, -n(K_X+\Delta+G)+f^*L)\to H^0(G, -n(K_G+\Delta|_G)+(f^*L)|_G). $$
Thus there is an element
$$nF\in |-n(K_X+\Delta+f^*A)+f^*L|$$ 
such that $nF|_G=nF_G$. Since $(G, F_G+\Delta|_G)$ is lc, by inversion of adjunction, we see  $(X, F+\Delta)$ is lc over a neighborhood of $H$. It is clear by construction that 
$$n(K_X+F+\Delta)\sim 0/Z.$$
Recall that $H$ is general in the linear system $|A|$. Taking a general element 
$$n\Theta\in |-n(K_X+\Delta+f^*A)+f^*L|,$$
we see that there is an open subset $U\subset Z$ with $\dim(Z\setminus U)=0$ satisfying that $(X, \Theta+\Delta)$ is lc over $U$. It is clear that $$n(K_X+\Theta+\Delta)\sim 0/Z.$$

\

\textit{Step 2}. Write $Z\setminus U:=\{z_1,...,z_k\}$, where $z_1,...,z_k$ are closed points on $Z$. By the construction in Step 1, we know that $(X, \Theta+\Delta)$ loses lc singularities on the fibers over these $z_i$. For each $z_i$, we denote by $V_i$ the sub-linear system $|f^*2A|$ consisting of elements containing the fiber $f^*{z_i}$. Since $A$ is very ample, we see $\dim V_i>0$. Let $p$ be a natural number such that $\frac{1}{p}<1-\epsilon$. Pick distinct general elements $M_{i,1}, M_{i,2},...,M_{i,p(d+1)}$ in $V_i$ and let 
$$M_i:=\frac{1}{p}(M_{i,1}+M_{i,2}+...+M_{i, p(d+1)}). $$
Then $(X, M_i+\Delta)$ is $\epsilon$-lc outside $f^*{z_i}$ by the generality of $M_{i,j}, 1\leq j\leq p(d+1)$. On the other hand, $(X, M_i+\Delta)$ is not lc at any point of $f^*z_i$ by \cite[Theorem 18.22]{Kollar92}.

Let $\epsilon'\in (0,\epsilon)$ be a sufficient small rational number and let $u$ be the largest number such that $(X, uM_i+\Delta)$ is $\epsilon'$-lc. Thus there is a prime divisor $E_i$ over $X$ such that there exists a birational contraction $g_i: Y_i\to X$ which extracts $E_i$ and
$$K_{Y_i}+u\widetilde{M_i}+\widetilde{\Delta}+(1-\epsilon')E_i=g_i^*(K_X+uM_i+\Delta). $$
Note that $Y_i$ is still Fano type over $Z$. We may assume $E_i$ is $g_i$-exceptional since the case when $E_i$ is a prime divisor on $X$ can be done similarly below.
Up to a small $\bQ$-factorialization of $Y_i$, we may assume $Y_i$ is $\bQ$-factorial. Observing that
$$-K_{Y_i}-\widetilde{\Delta} -E_i =-g_i^*(K_X+uM_i+\Delta)+u\widetilde{M_i}-\epsilon'E_i$$
is big over $Z$ (since $E_i$ is vertical over $Z$), we run MMP/$Z$ on $-K_{Y_i}-\widetilde{\Delta}-E_i $ to get a minimal model $Y_i\dashrightarrow Y_i'/Z$ such that $-K_{Y_i'}-\widetilde{\Delta}'-E_i' $ is big and nef over $Z$. We claim that $E_i$ is not contracted by $Y_i\dashrightarrow Y_i'$. Indeed, let $$0<P\sim_\bQ -(K_X+\Delta)/ Z$$ 
be a general ample $\bQ$ divisor such that $(X, uM_i+\Delta+P)$ is still $\epsilon'$-lc and we have
$$K_{Y_i}+u\widetilde{M_i}+\widetilde{\Delta}+g_i^*P+(1-\epsilon')E_i=g_i^*(K_X+uM_i+\Delta+P). $$
Then $(Y_i,u\widetilde{M_i}+\widetilde{\Delta}+g_i^*P+(1-\epsilon')E_i)$ is also $\epsilon'$-lc. Since 
$$K_X+uM_i+\Delta+P\sim_\bQ 0/Z,$$ 
the log pair 
$$(Y'_i,u\widetilde{M_i}'+\widetilde{\Delta}'+(g_i^*P)'+(1-\epsilon')E'_i)$$ 
is $\epsilon'$-lc and so is $(Y_i', \widetilde{\Delta}'+(1-\epsilon')E_i')$. Suppose $E_i$ is contracted. Recall that $Y_i\dashrightarrow Y_i'$ is obtained by MMP/$Z$ on $-K_{Y_i}-\widetilde{\Delta}-E_i$. Thus we see
$$0=A_{(Y_i, \widetilde{\Delta}+E_i)}(E_i)>A_{(Y_i', \widetilde{\Delta}')}(E_i)\geq \epsilon', $$
which leads to a contradiction. Hence $E_i$ is not contracted.

We may choose $\epsilon'$ sufficiently small such that $(Y'_i, E'_i)$ is lc (e.g. \cite{HMX14}). 
Write $f_i': Y_i'\to Z$. 
Let $A_i$ be a very ample Cartier divisor on $Z$ such that
$$f_i'^*A_i -K_{Y_i'}-\widetilde{\Delta}'-E_i' $$
is big and nef. Now we can apply Theorem \ref{thm: global engine} to see there exist a positive integer $m$ (depending only on $d, \epsilon, \mR$) and an effective $\bQ$-divisor $\Lambda_i\geq  \widetilde{\Delta}'+E_i'$ satisfying 
$$m(K_{Y_i'}+\Lambda_i)\sim (m+2)f_i'^*A_i$$ 
and $(Y_i', \Lambda_i)$ is lc over $z_i$. By \cite[6.1(3)]{Birkar19}, there exists an effective $\bQ$-divisor $\Omega_i\geq \widetilde{\Delta}+E_i$ such that 
$$m(K_{Y_i}+\Omega_i)\sim (m+2)g_i^*f^*A_i$$ 
and $(Y_i, \Omega_i)$ is lc over $z_i$. Denote by $\Theta_i:={g_i}_*\Omega_i\geq \Delta$. It is clear that 
$$m(K_X+\Theta_i)\sim (m+2)f^*A_i$$ 
and $(X, \Theta_i)$ is lc over $z_i$.

\

\textit{Step 3}. Up to replacing $m$ and $m+2$ by adding a common multiple, we may assume that $m$ is a multiple of $n$ (see Step 1 for $n$) and write $m=rn$. By now we have
\begin{enumerate}
\item $m(\Theta_i-\Delta)\in |(m+2)f^*A_i-m(K_X+\Delta)| ,$
\item $(X, \Theta_i)=(X, (\Theta_i-\Delta)+\Delta)$ is lc over $z_i$,
\item $n\Theta\in |-n(K_X+\Delta)-nf^*A+f^*L|$,
\item $(X, \Theta+\Delta)$ is lc over $X\setminus \{z_1,...,z_k\}$.
\end{enumerate}
Note that the last two items are obtained in Step 1. Let $M$ be a very ample Cartier divisor on $Z$ such that 
$$M-(m+2)A_i $$
is very ample for every $i$, and 
$$M+rnA-rL $$
is very ample. Then we see
$$m(\Theta_i-\Delta)+f^*M-(m+2)f^*A_i\in |f^*M-m(K_X+\Delta)| $$
and
$$m\Theta+f^*M+rnf^*A-rf^*L\in |f^*M-m(K_X+\Delta)|.  $$
Let $m\Lambda\in |f^*M-m(K_X+\Delta)|$ be a general element in the linear system, then we have
\begin{enumerate}
\item $(X, \Lambda+\Delta)$ is lc over $z_i$ for each $i$,
\item $(X, \Lambda+\Delta)$ is lc over $X\setminus \{z_1,...,z_k\}$.
\end{enumerate}
Thus the two items above together imply that $(X, \Lambda+\Delta)$ is lc. It is clear that 
$$m(K_X+\Lambda+\Delta)\sim 0/Z,$$ 
and $\Lambda+\Delta$ is the required bounded relative-global complement.
\end{proof}

\begin{remark}\label{rem: what we lack}
One may notice that we only partially apply Theorem \ref{thm: curve} for induction in the proof of Theorem \ref{thm: main}. That is to say, we prove the existence of bounded relative-global klt complements for curve base in Theorem \ref{thm: curve} while we only apply the existence of bounded relative-global (lc) complements for induction. If we fully apply Theorem \ref{thm: curve} for induction, the same proof of Theorem \ref{thm: main} produces a bounded relative-global complement which is klt outside finitely many closed fibers. This is quite close to the final solution to Conjecture \ref{conj: main}, but it seems difficult to construct a bounded relative-global complement which is klt over the finitely many closed points.
\end{remark}

\section{Final remark}

In this section, we prove Theorem \ref{main: global lc} for Fano fibrations  by a more direct way without using induction.

\begin{theorem}
Let $d$ be a positive integer and $\mR\subset[0,1]$ be a finite set of rational numbers. Then there exists a positive integer $N$ depending only on $d, \mR$ satisfying the following. 
Let $f: (X, \Delta)\to Z$ be a log Fano fibration of dimension $d$ where the coefficients of $\Delta$ are contained in $\mR$. Then there exists an effective $\bQ$-divisor $\Lambda\geq \Delta$ on $X$ satisfying that $N(K_X+\Lambda)\sim 0/Z$ and $(X, \Lambda)$ is lc.
\end{theorem}

\begin{proof}
Fix a closed point $z\in Z$. By the construction in Step 2 of the proof of Theorem \ref{thm: main},  there exist a positive integer $m$ (depending only on $d, \mR$) and an effective $\bQ$-divisor $\Theta\geq \Delta$ such that 
$$m(K_X+\Theta)\sim 0/Z$$ 
and $(X, \Theta)$ is lc over $z$. We may choose $m$ sufficiently divisible such that $m\Delta$ is integral.

Given this, we may find finitely many closed points $z_1,...,z_k$ on $Z$, and for each $z_i$ there is an effective $\bQ$-divisor $\Theta_i\geq \Delta$ satisfying 
$$m(K_X+\Theta_i)\sim 0/Z$$ 
and $(X, \Theta_i)$ is lc over a neighborhood $U_i$ of $z_i$. We may assume $\cup_i U_i=Z$ and write 
$$m(K_X+\Theta_i)\sim f^* A_i$$ 
for some Cartier divisor $A_i$ on $Z$. Then it is clear that 
$$m(\Theta_i-\Delta)\in |f^*A_i-m(K_X+\Delta)|.$$

Let $M$ be a very ample Cartier divisor on $Z$ such that $M-A_i$
is very ample for every $i$. Then we see
$$m(\Theta_i-\Delta)+f^*M-f^*A_i\in |f^*M-m(K_X+\Delta)|. $$
Let $m\Lambda\in |f^*M-m(K_X+\Delta)|$ be a general element in the linear system. Then we see $(X, \Lambda+\Delta)$ is lc over $U_i$ for each $i$.
This implies that $(X, \Lambda+\Delta)$ is lc. It is clear that 
$$m(K_X+\Lambda+\Delta)\sim 0/Z,$$ 
and $\Lambda+\Delta$ is the required bounded relative-global complement.
\end{proof}

We end this note by the following remark.

\begin{remark}
All the results presented here for Fano fibrations (resp. $\epsilon$-lc Fano fibrations) can be generalized to Fano type fibrations (resp. $\epsilon$-lc Fano type fibrations) without any difficulty. For example, let $X\to Z$ be a Fano type fibration (resp. $\epsilon$-lc Fano type fibration), then there exists an effective $\bQ$-divisor $B$ on $X$ such that $(X, B)\to Z$ is a log Fano fibration (resp. $\epsilon$-lc log Fano fibration). In particular, $-K_X$ is big over $Z$. Running MMP/$Z$ on $-K_X$ to get a good minimal model $X\dashrightarrow X'/Z$ and then considering the ample model $X'\dashrightarrow X''/Z$ where $-K_{X''}$ is ample over $Z$, it is not hard to see that $X''\to Z$ is a Fano fibration (resp. $\epsilon$-lc Fano fibration). By \cite[6.1(3)]{Birkar19}, we see all the results on complements established for $X''\to Z$ also hold for $X\to Z$.
\end{remark}

\bibliography{reference.bib}
\end{document}